\documentclass[12pt]{article}
  \usepackage[cp1251]{inputenc}
  \usepackage[english]{babel}
  \usepackage{amsfonts,amssymb,amsmath,amstext,amsbsy,amsopn,amscd,graphicx}
  \usepackage{graphics}
  \usepackage{amsthm}

  \def\thtext#1{
  \catcode`@=11
  \gdef\@thmcountersep{. #1}
  \catcode`@=12}

  \newtheorem{theorem}{Theorem}[section]
  
  \newtheorem{st}{Statement}[section]

 \newcounter{il}[section]
 \renewcommand{\thtext}
 {\thechapter.\arabic{il}}

 \newcounter{il2}[section]
 \renewcommand{\thtext}
 {\thechapter.\arabic{il2}}

 \newenvironment{dfn}{\trivlist \item[\hskip\labelsep{\bf Definition}]
 \refstepcounter{il}{\bf \arabic{section}.\arabic{il}.}}%
 {\endtrivlist}

 \newenvironment{rk}{\trivlist \item[\hskip\labelsep{\bf Remark}]
 \refstepcounter{il2}{\bf \arabic{section}.\arabic{il2}.}}%
 {\endtrivlist}

 \newenvironment{examp}{\trivlist \item[\hskip\labelsep{\bf Example.}]}%
 {\endtrivlist}

 {\endtrivlist}

 \catcode`@=11
 \def\.{.\spacefactor\@m}
\def\relaxnext@{\let\next\relax}
\def\nolimits@{\relaxnext@
 \DN@{\ifx\next\limits\DN@\limits{\nolimits}\else
  \let\next@\nolimits\fi\next@}%
 \FN@\next@}
\def\newmcodes@{\mathcode`\'39\mathcode`\*42\mathcode`\."613A%
 \mathcode`\-45\mathcode`\/47\mathcode`\:"603A\relax}
\def\operatorname#1{\mathop{\newmcodes@\kern\z@
 \operator@font#1}\nolimits@}
 \catcode`@=12
 \def\rom#1{{\em#1}}

 \def\({\rom(}
 \def\){\rom)}
 \def\:{\colon}

 \def\Z{{\mathbb Z}}
 \def\0{{\mathbf 0}}
 \def\1{{\mathbf 1}}

\newcommand{\chord}{\raisebox{-0.25\height}{\includegraphics[width=0.7cm]{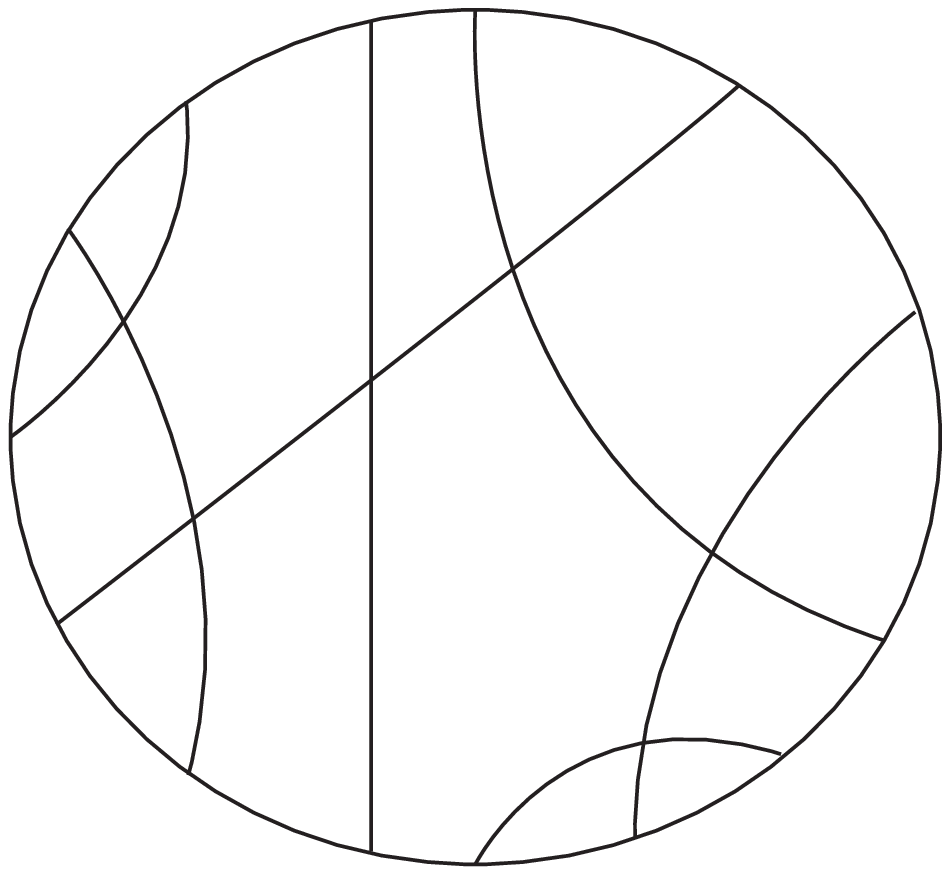}}}

\newcommand{\chordff}{\raisebox{-0.25\height}{\includegraphics[width=1.65cm]{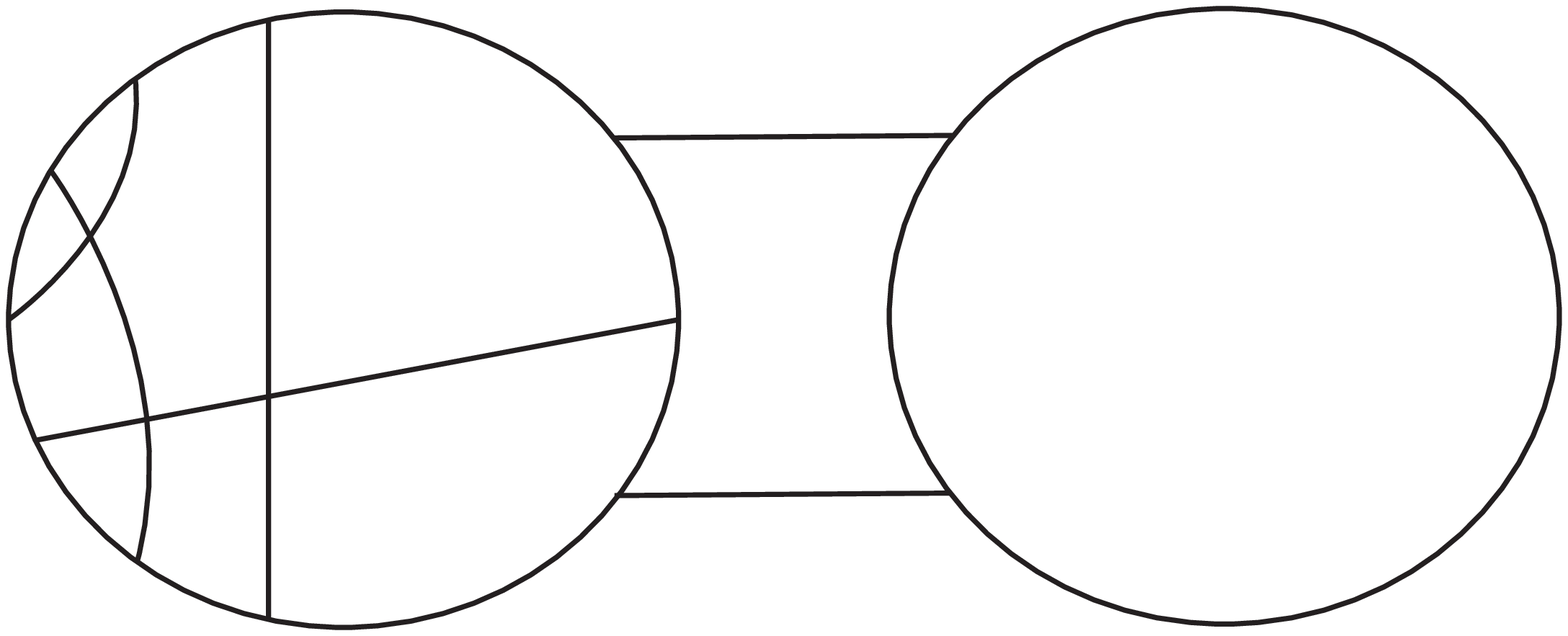}}}
\newcommand{\chordfs}{\raisebox{-0.25\height}{\includegraphics[width=1.65cm]{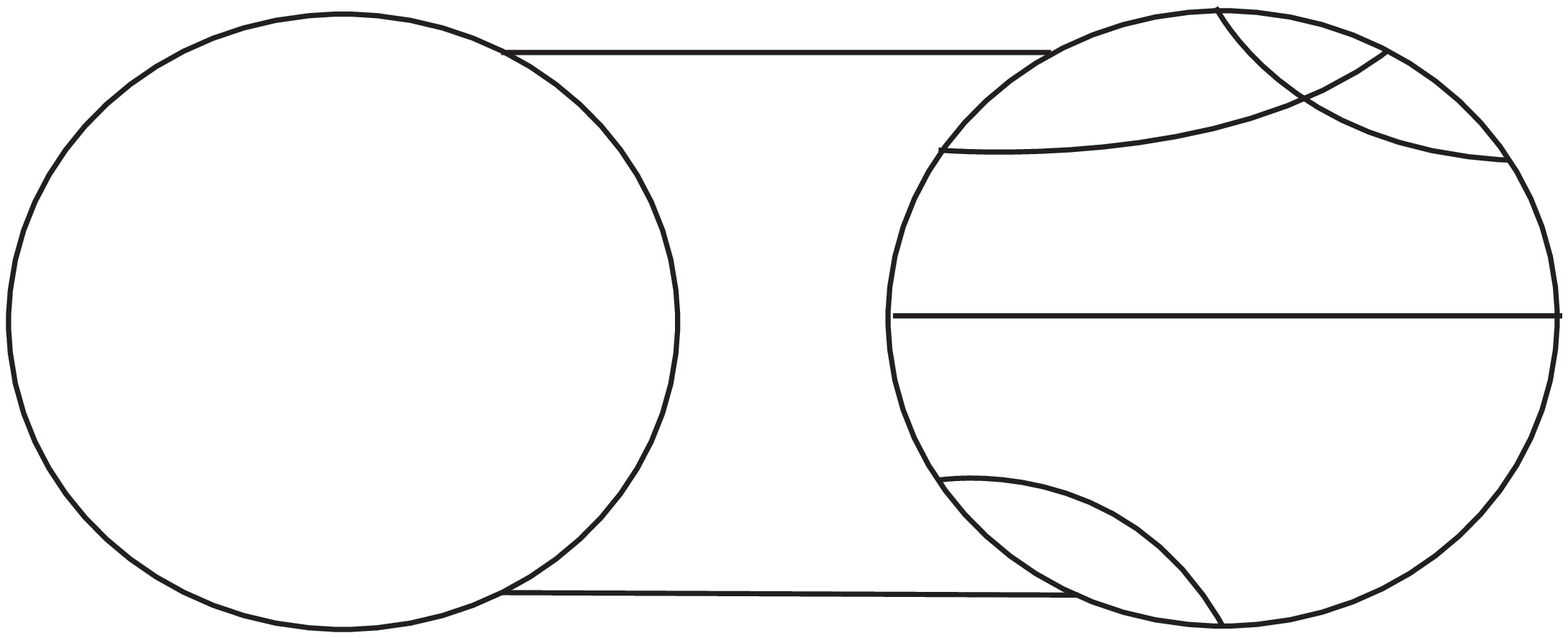}}}
\newcommand{\chordft}{\raisebox{-0.25\height}{\includegraphics[width=1.65cm]{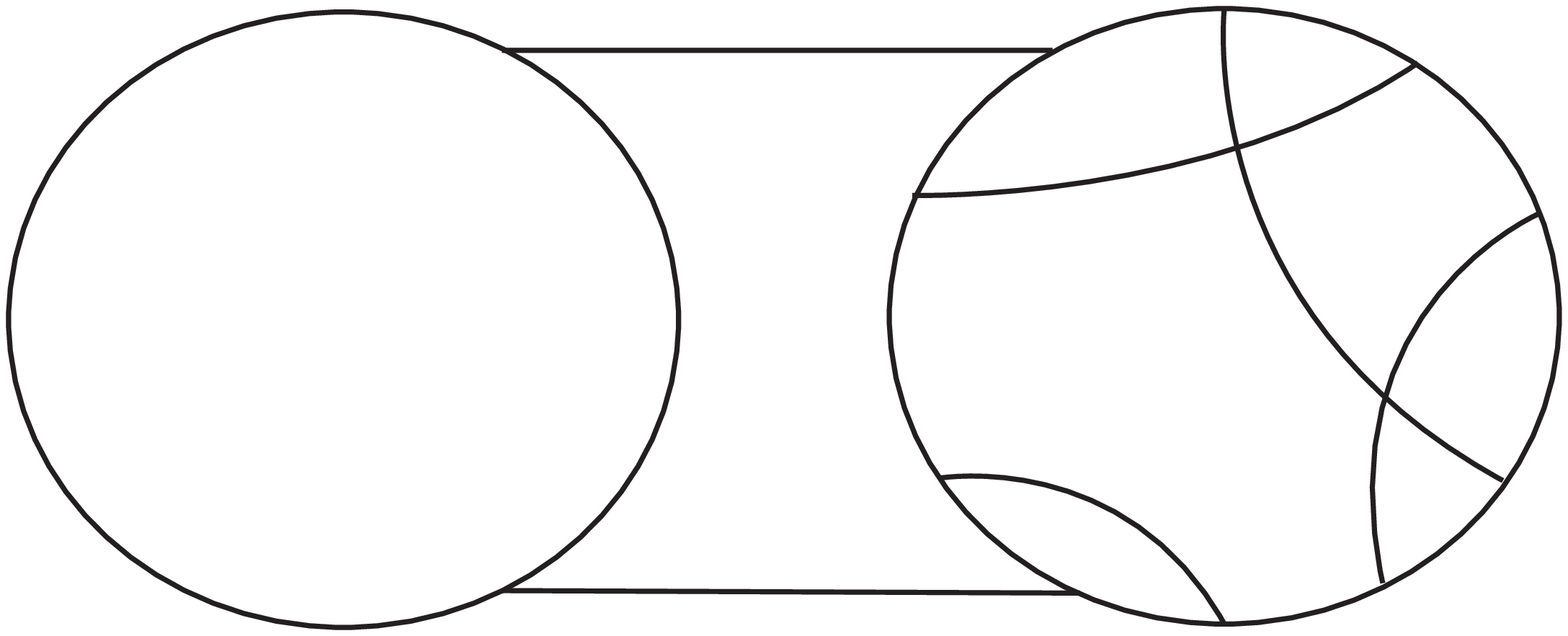}}}

\newcommand{\chordsf}{\raisebox{-0.4\height}{\includegraphics[width=1.65cm]{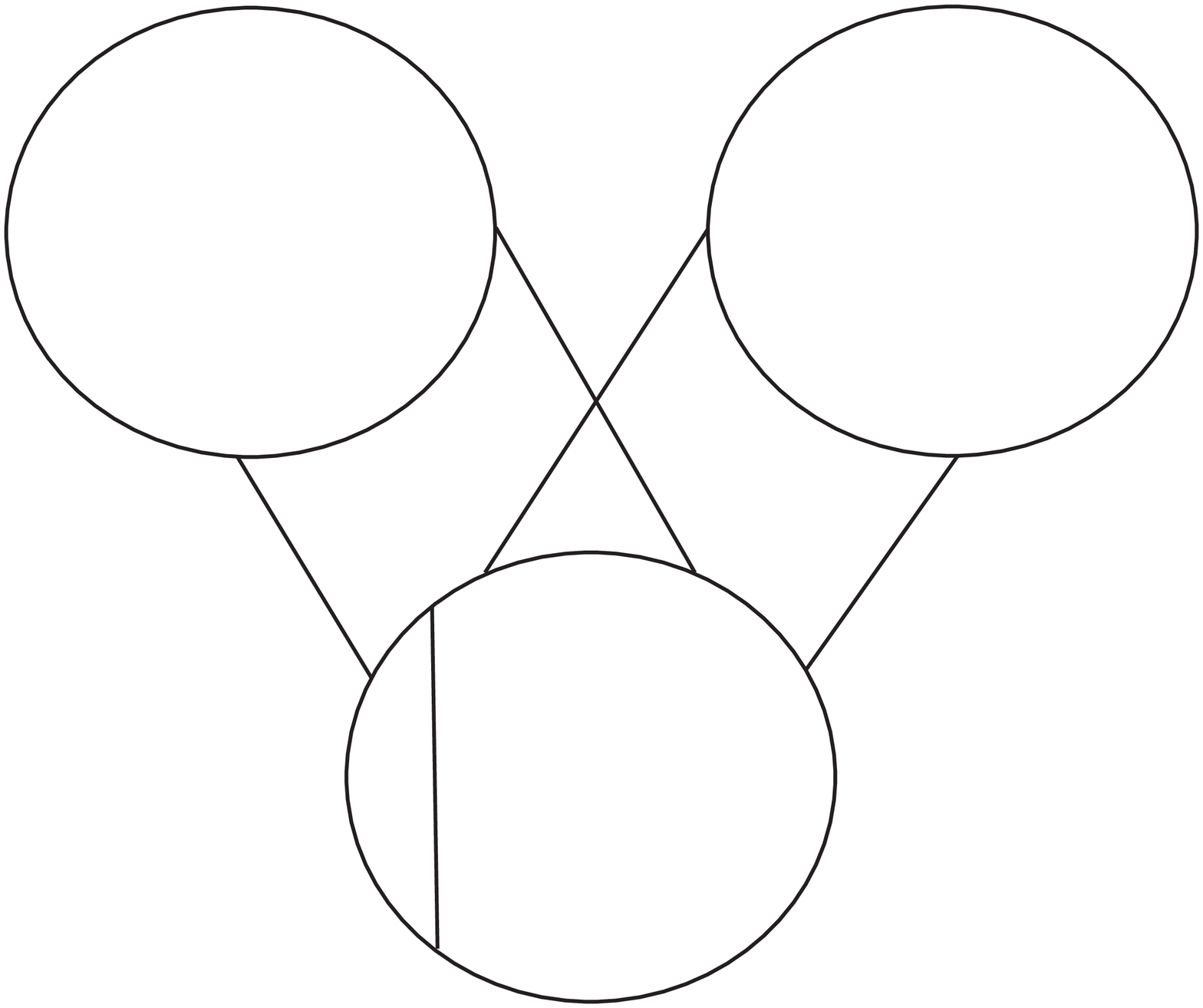}}}
\newcommand{\chordss}{\raisebox{-0.4\height}{\includegraphics[width=1.65cm]{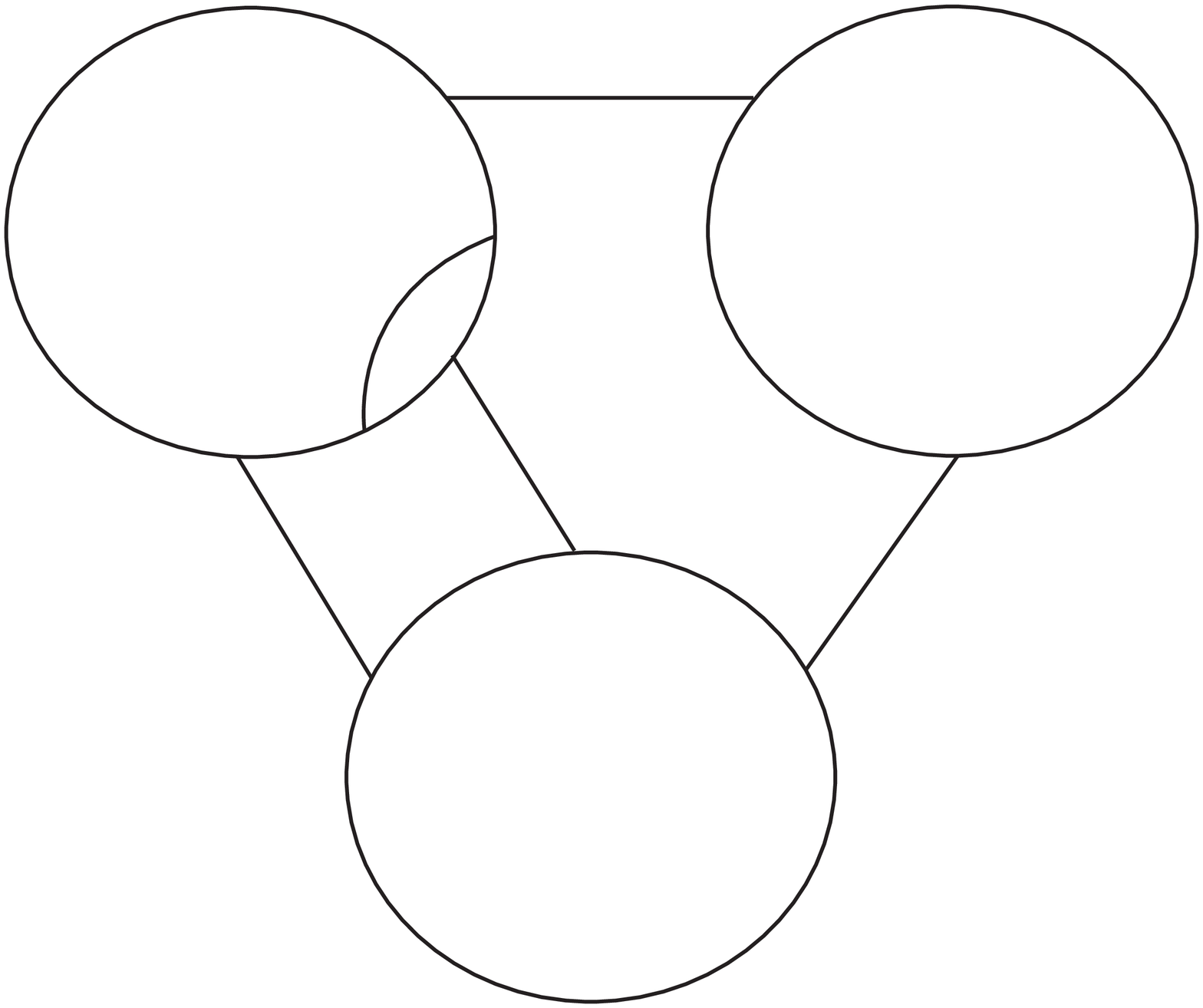}}}
\newcommand{\chordst}{\raisebox{-0.4\height}{\includegraphics[width=1.65cm]{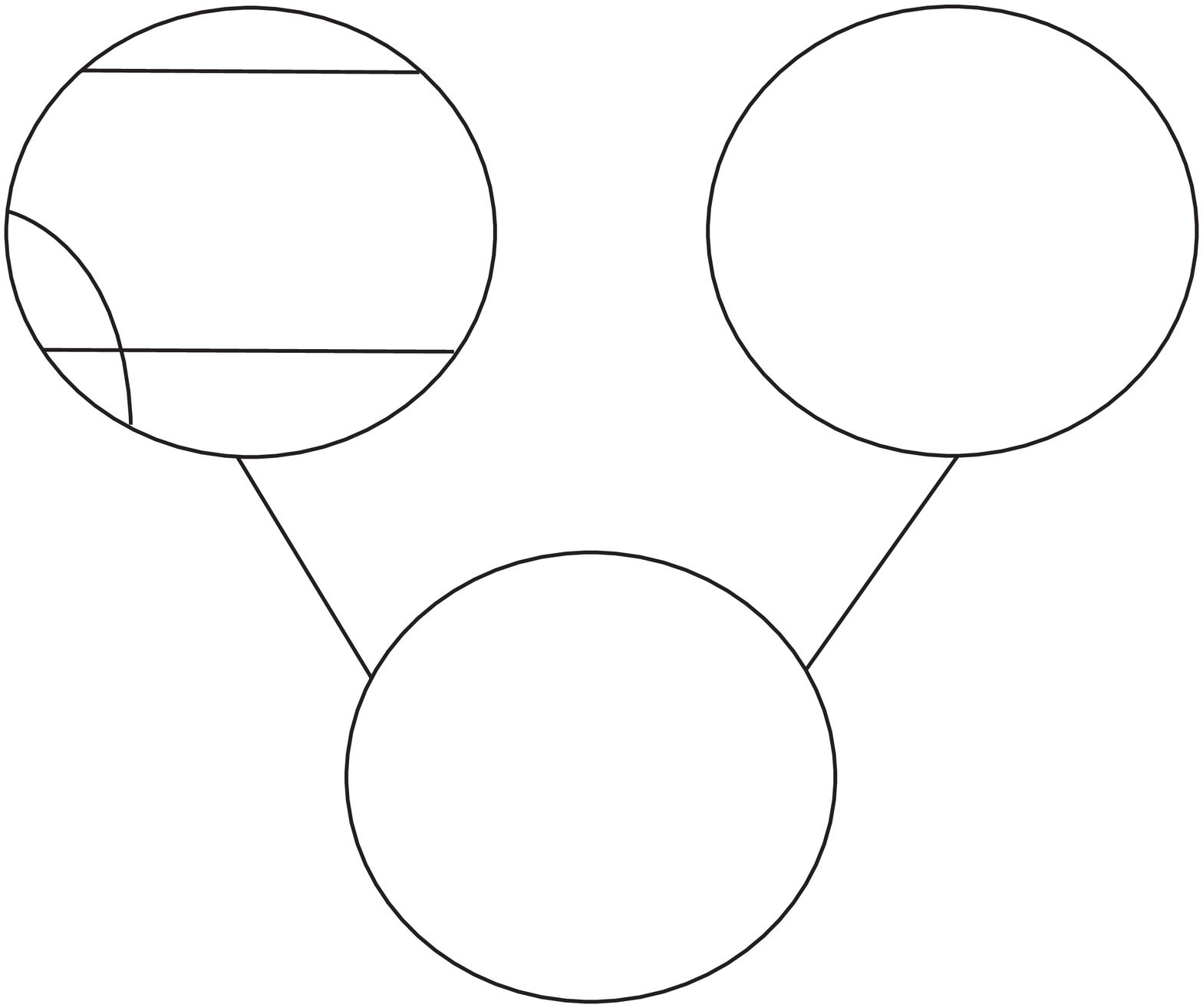}}}
\newcommand{\chordsfo}{\raisebox{-0.4\height}{\includegraphics[width=1.65cm]{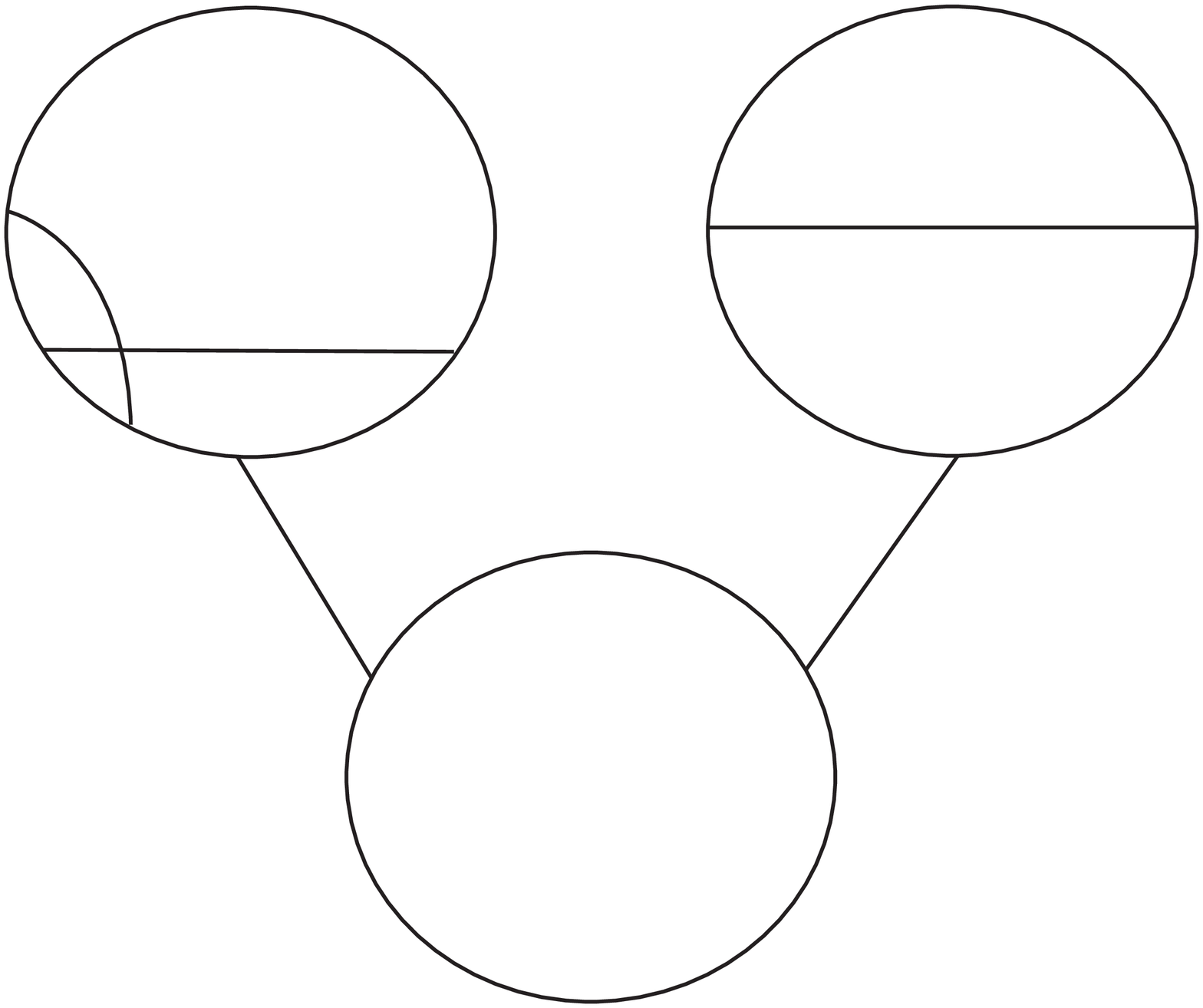}}}
\newcommand{\chordsfi}{\raisebox{-0.4\height}{\includegraphics[width=1.65cm]{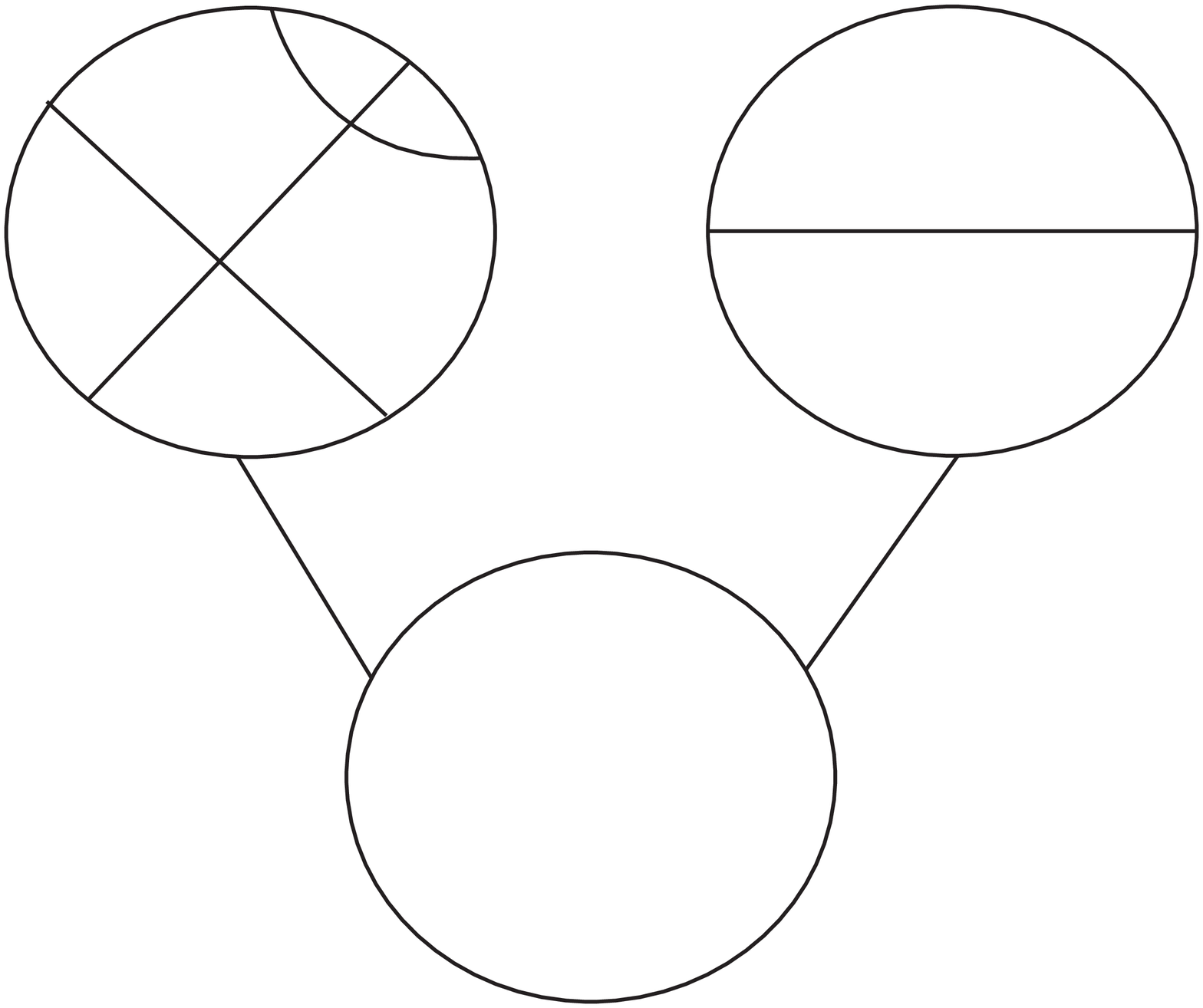}}}
\newcommand{\chordssi}{\raisebox{-0.4\height}{\includegraphics[width=1.65cm]{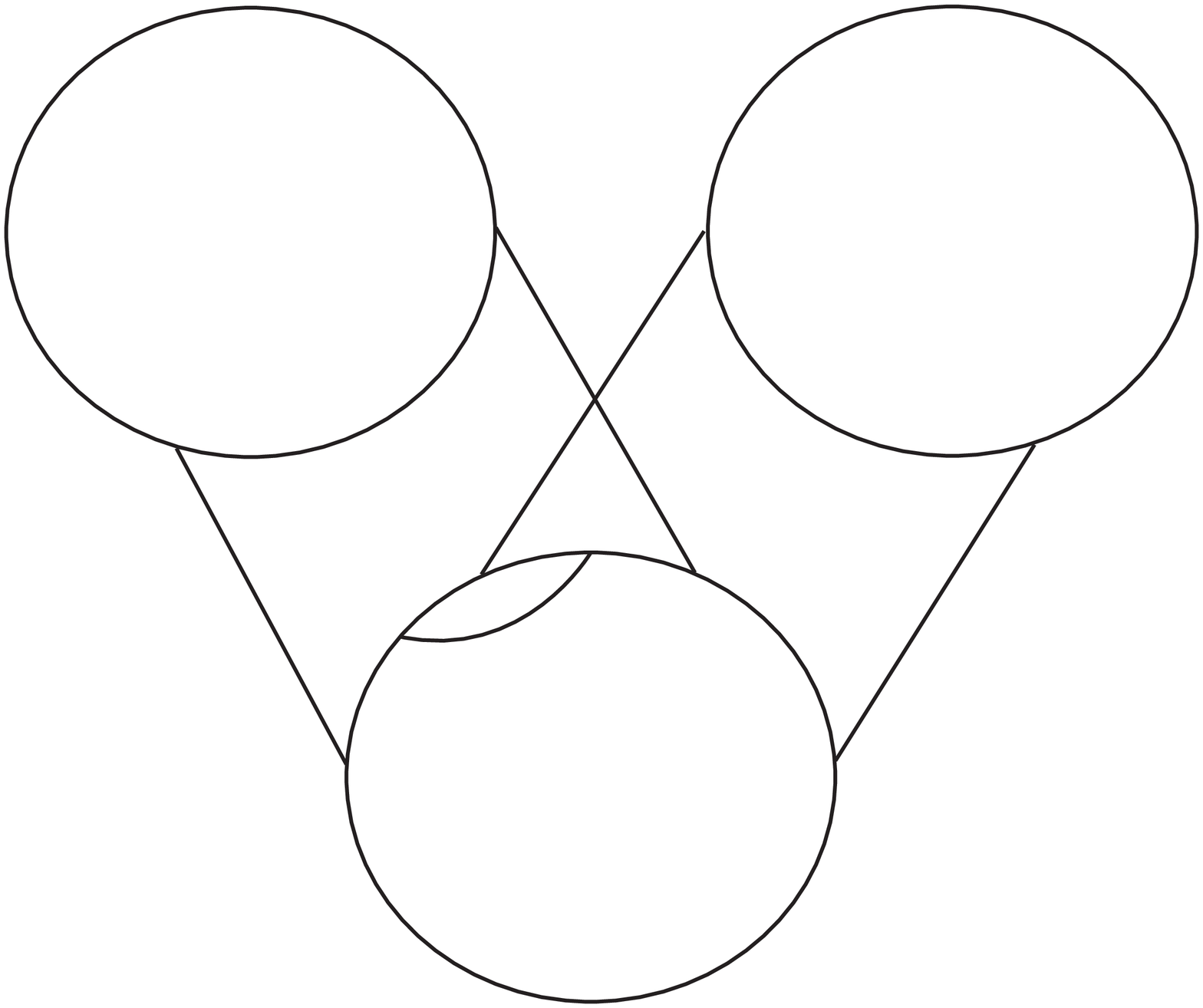}}}

\newcommand{\chordtf}{\raisebox{-0.4\height}{\includegraphics[width=1.65cm]{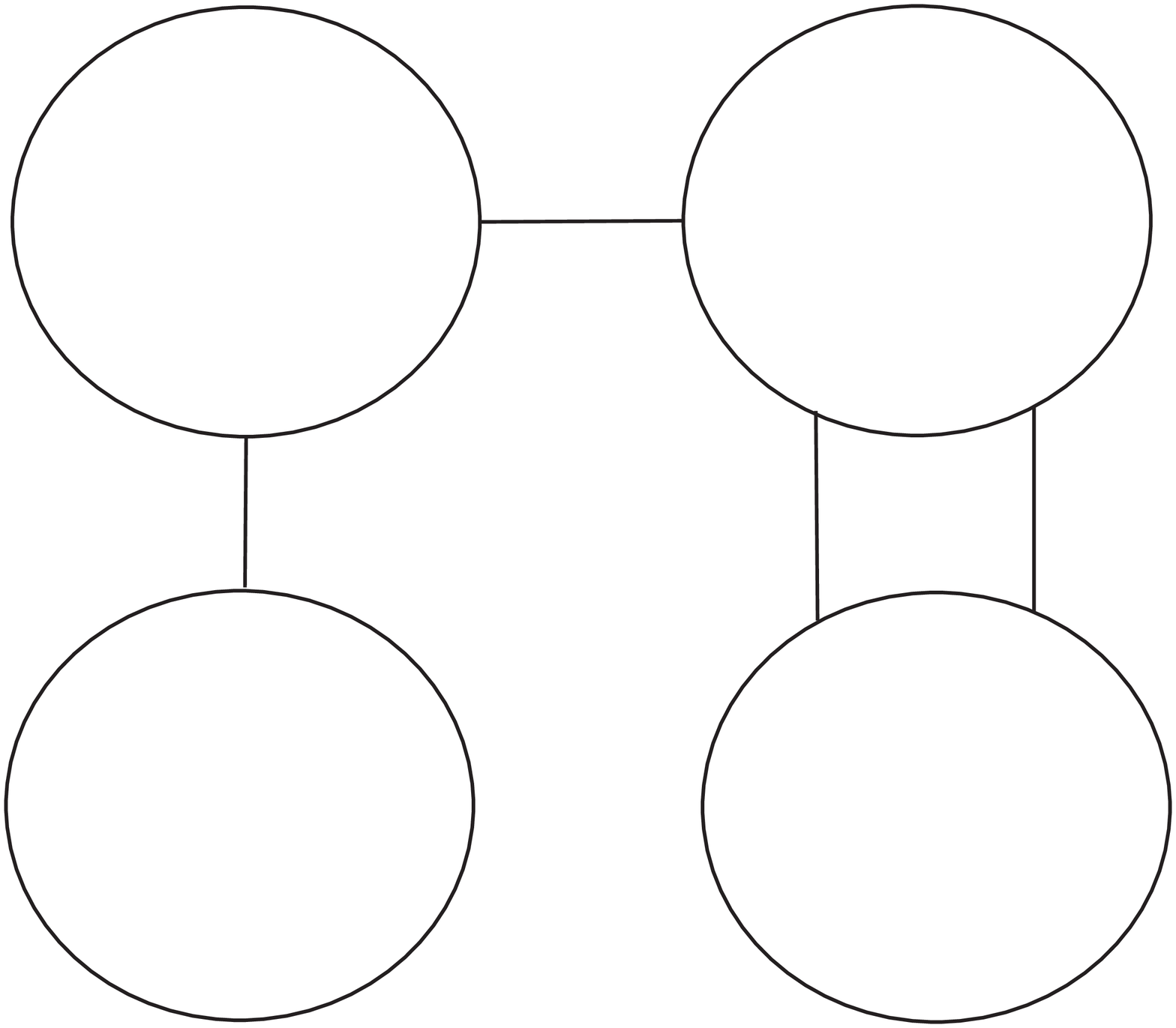}}}
\newcommand{\chordts}{\raisebox{-0.4\height}{\includegraphics[width=1.65cm]{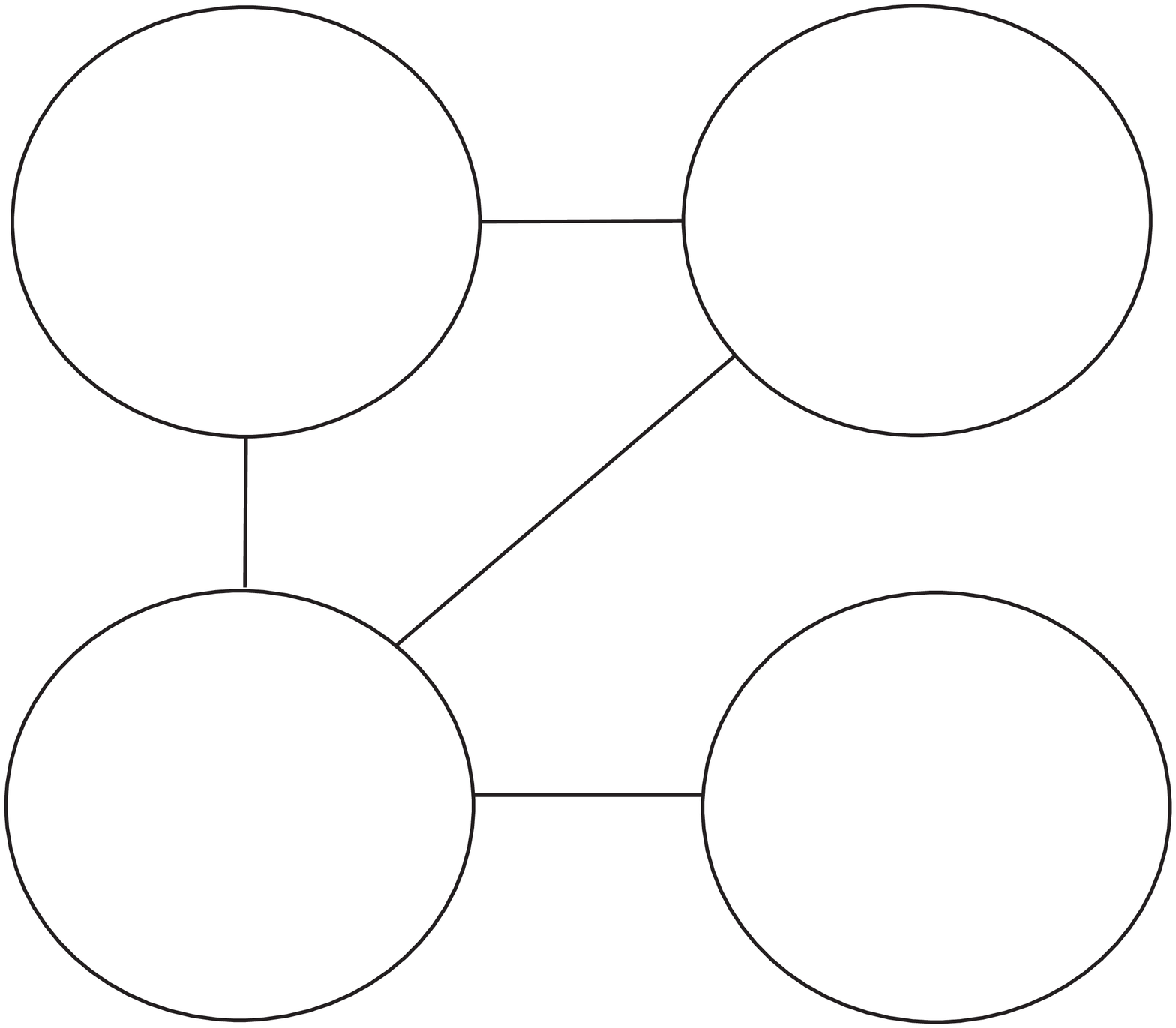}}}
\newcommand{\chordtt}{\raisebox{-0.4\height}{\includegraphics[width=1.65cm]{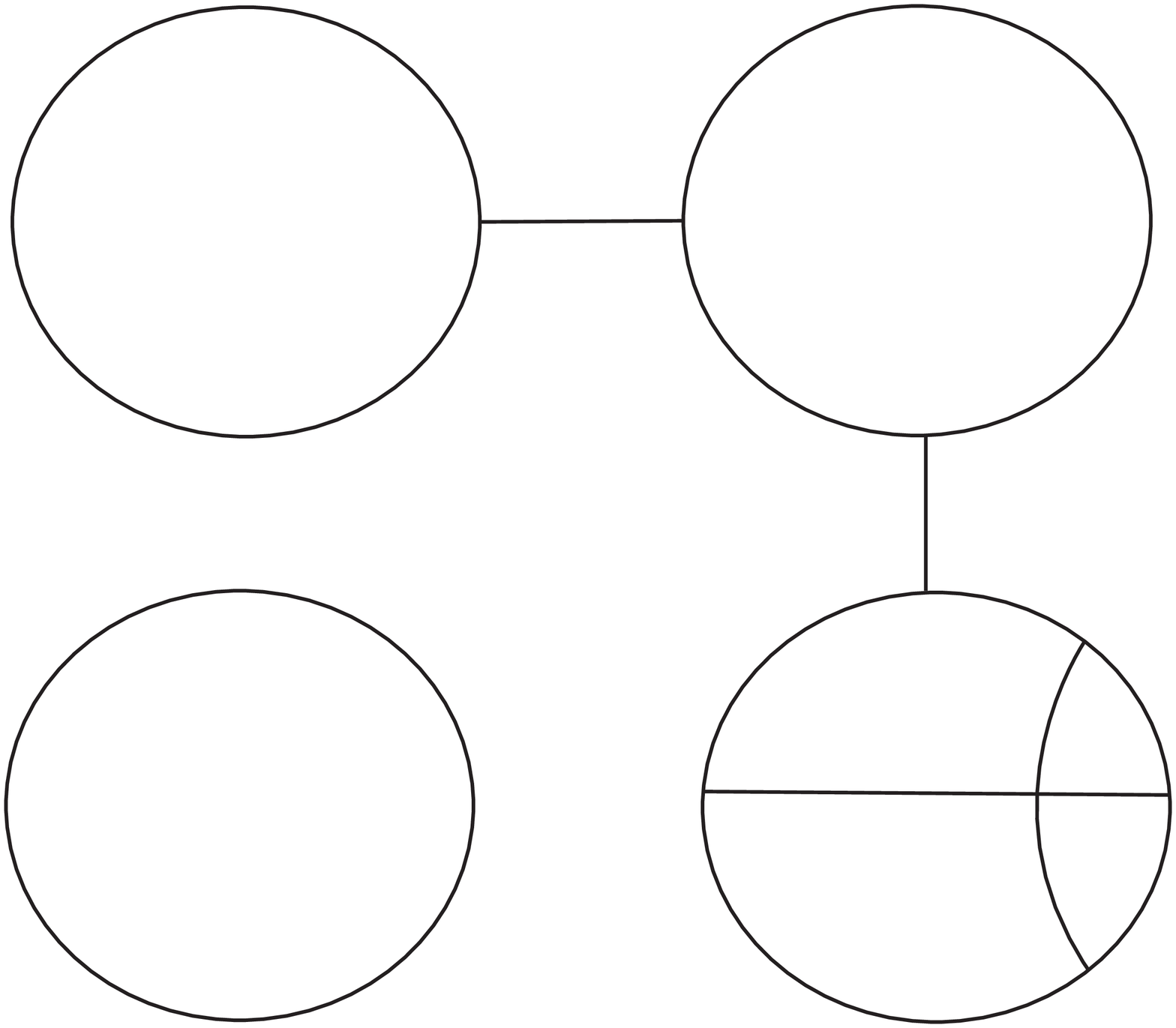}}}
\newcommand{\chordtfo}{\raisebox{-0.4\height}{\includegraphics[width=1.65cm]{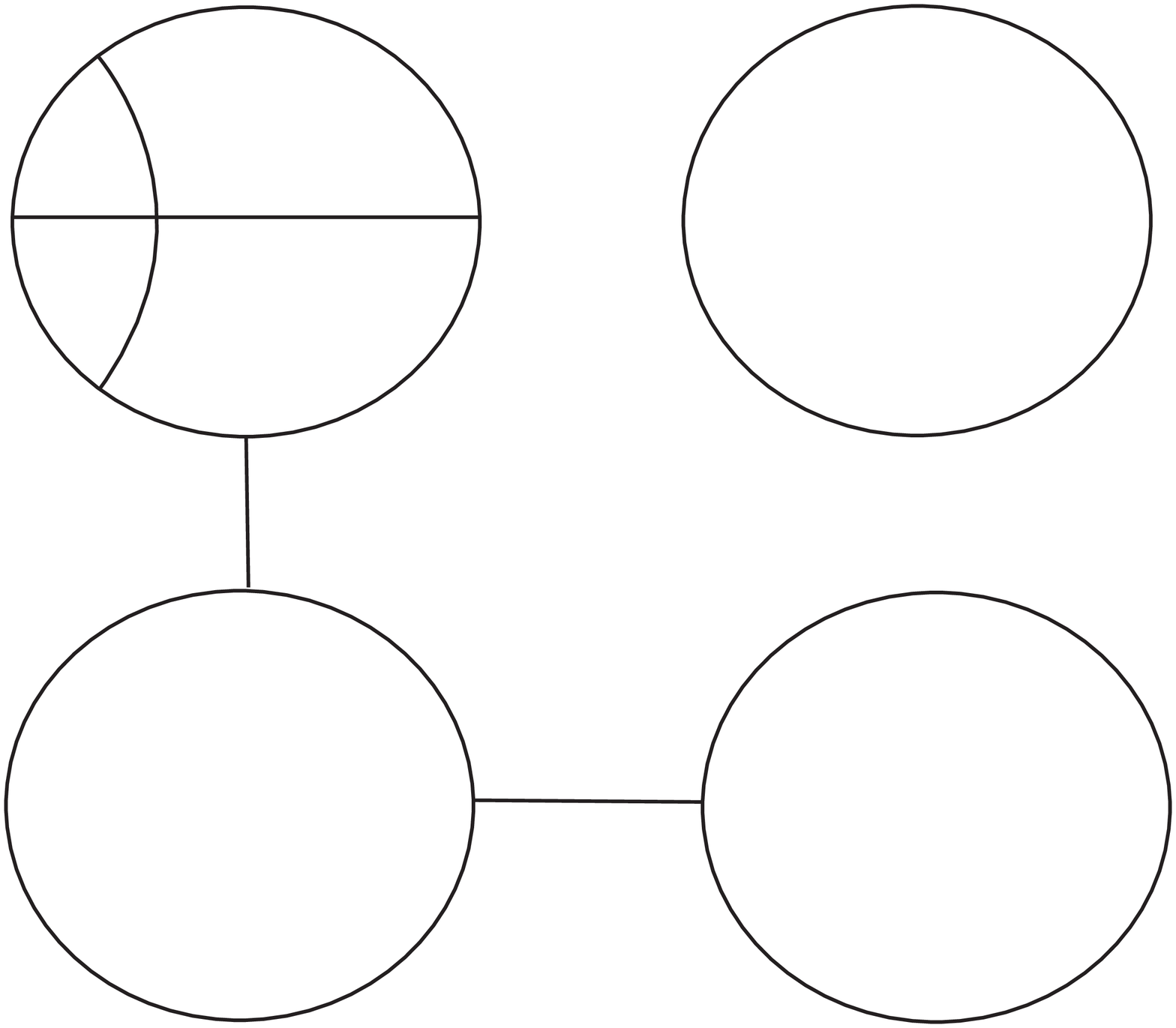}}}
\newcommand{\chordtfi}{\raisebox{-0.4\height}{\includegraphics[width=1.65cm]{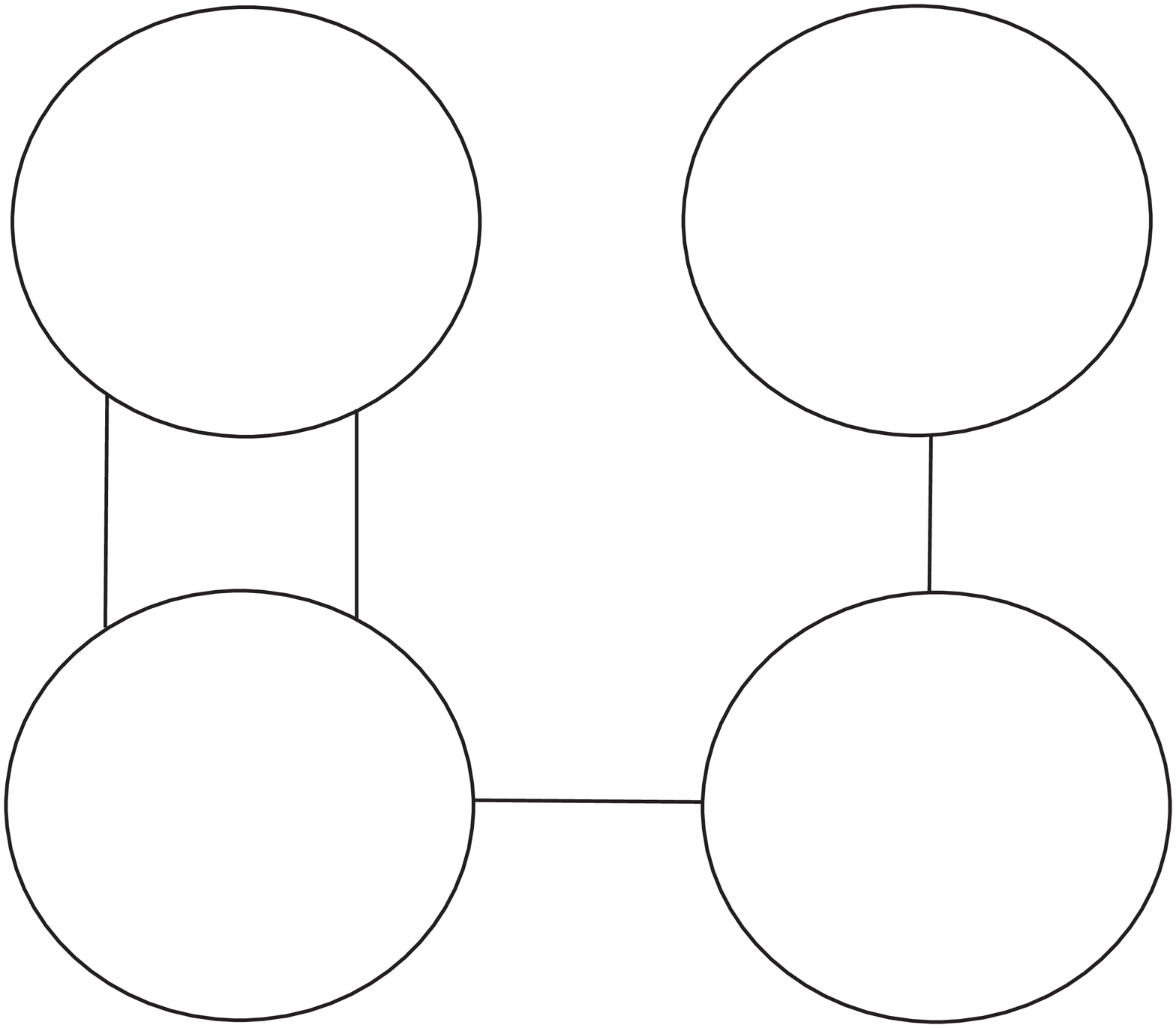}}}
\newcommand{\chordtsi}{\raisebox{-0.4\height}{\includegraphics[width=1.65cm]{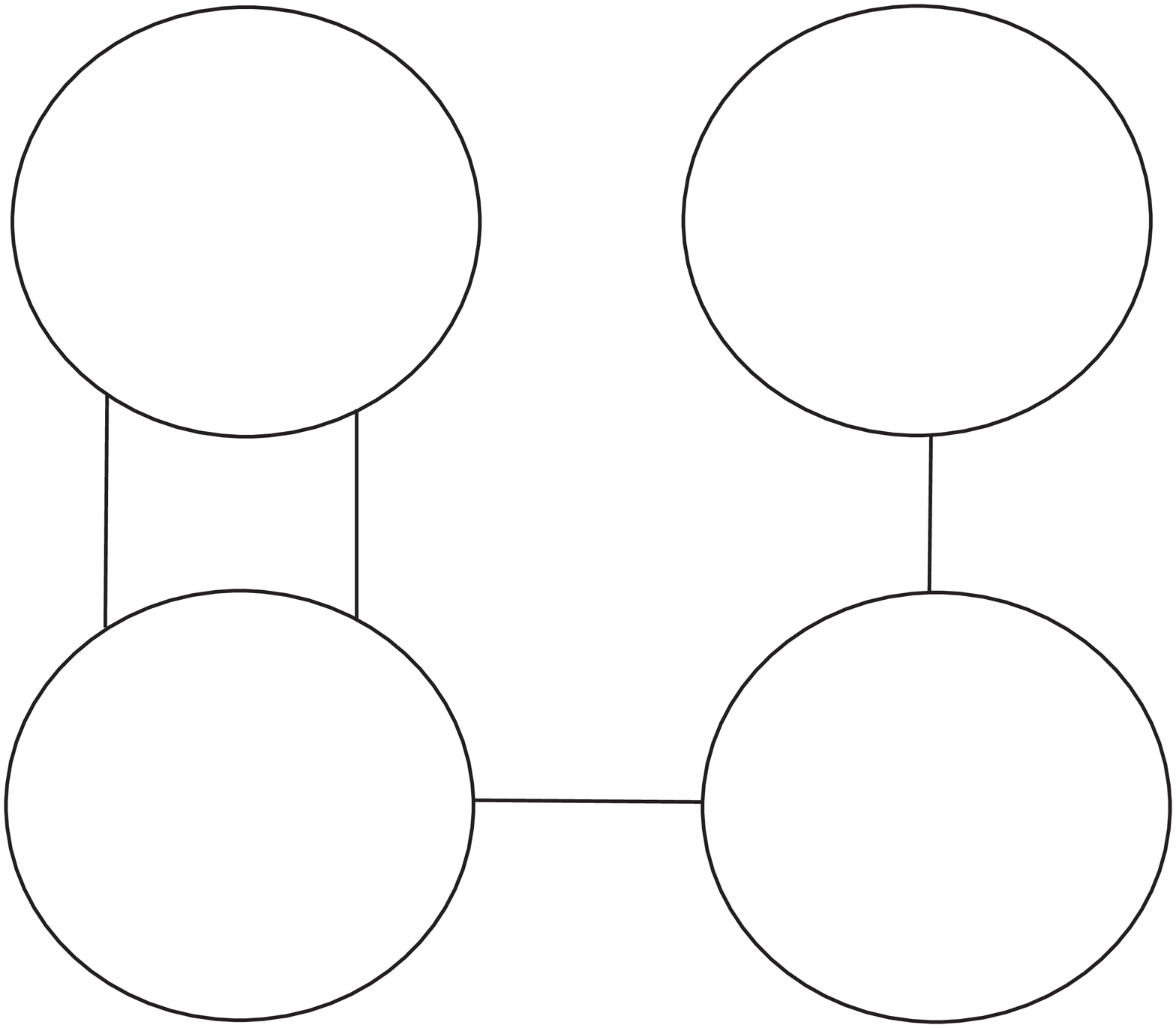}}}

 \title {Cobordisms of Free Knots and Gauss Words}

\author{D.\,P.~Ilyutko\footnote{Partially supported by grants of RF President
NSh -- 660.2008.1, RFFI 07--01--00648, RNP 2.1.1.7988.}\,,
V.\,O.~Manturov}
\date {}

\begin {document}

 \maketitle
 \abstract{We investigate cobordisms of free knots. Free knots and links are also called homotopy
classes of Gauss words and phrases. We define a new strong invariant
of free knots which allows to detect free knots not cobordant to the
trivial one.}


 \section {Introduction}


The aim of the present work is to consider {\em cobordisms of free
knots}. Free knots and links \cite{Ma} (also called {\em homotopy
classes of Gauss words and phrases}, see \cite{Tu,Gib}) are a
substantial simplification of homotopy classes of curves on
$2$-surfaces, and at the same time, a simplification of virtual
knots and links, introduced by Kauffman in \cite{Kau}.

A conjecture due to Turaev \cite{Tu} stating that all free knots are
trivial was solved very recently  independently by Manturov
\cite{Ma} and Gibson \cite{Gib}.

We consider {\em cobordism classes of free knots} and give a strong
invariant which allows to detect free knots not cobordant to the
trivial one.

\section{Basic Definitions}

We shall encode free knots and links by using framed graphs.

Throughout the paper by a {\em four-valent graph} we mean the
following ge\-ne\-ra\-li\-za\-tion: a finite $1$-complex $\Gamma$
with each connected component being homeomorphic either to a circle
or to a four-valent graph; by a {\em vertex} of a four-valent graph
we mean only vertices of those components which are homeomorphic to
four-valent graphs, and by edges we mean both edges of four-valent
graphs and circular components (the latter will be called {\em
cyclic edges}).

We call a four-valent graph {\em framed} if for every vertex of it a
way of splitting of the four incident half-edges into two pairs is
indicated. Half-edges belonging to the same pair are to be called
{\em opposite}. We shall also use the term {\em opposite} with
respect to {\em edges} containing opposite half-edges.

By an {\em isomorphism} of framed $4$-graphs we always mean a
framing-preserving homeomorphism.

By a {\em unicursal component} of a framed four-valent graph we mean
either its connected component homeomorphic to the circle or an
equivalence class of edges of some non-circular component, where the
equivalence is generated by the opposite relation.

Let  $K$ be a framed four-valent graph, and let $X$ be a vertex of
$X$. Let $(a,c)$ and $(b,d)$ be the two pairs of opposite half-edges
of the graph $K$ at $X$. By smoothings of $K$ at $X$ we call two
framed four-valent graphs obtained from $K$ by deleting $X$ and
reconnecting the adjacent (non-opposite) edges: for one graph we
connect the pairs $(a,b)$ and $(c,d)$, and for the other graph we
connect $(a,d)$ and $(b,c)$.

Consider a chord diagram\footnote{Usually in literature one
considers oriented and non-oriented singular knots and corresponding
oriented and non-oriented chord diagrams. In the present paper, we
restrict ourselves to the non-orientable case. Many theorems and
constructions can be easily extended to the orientable case.} $C$.
Then the corresponding framed four-valent graph $G(C)$ with a unique
unicursal component is constructed as follows. With the chord
diagram having no chords one associates the graph  $G_{0}$
consisting of a unique circular component. Otherwise the edges of
the graph are in one-to-one correspondence with arcs of the chord
diagrams, and vertices are in one-to-one correspondence with chords.

Those arcs which are incident to the same chord end, correspond to
formally opposite half-edges.

\begin{figure}
\centering\includegraphics[width=200pt]{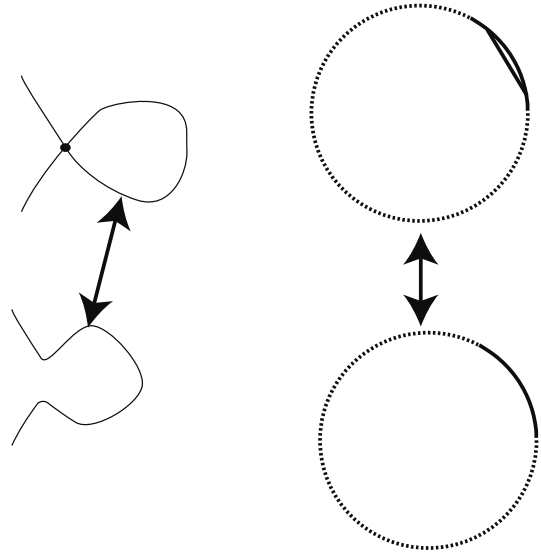} \caption{The First
Reidemeister Move and Its Chord Diagram Version} \label{1r}
\end{figure}

The first Reidemeister move for four-valent framed graphs is an
addition/re\-mo\-val of a loop, see Fig.~\ref{1r}.

\begin{figure}
\centering\includegraphics[width=200pt]{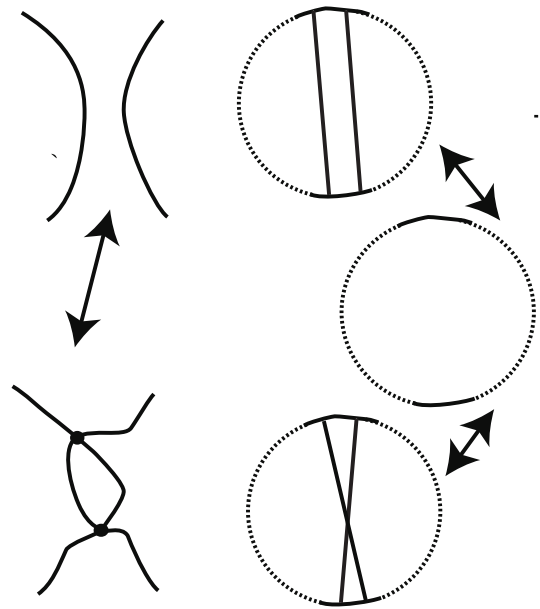} \caption{The Second
Reidemeister Move and Its Chord Diagram Version} \label{2r}
\end{figure}

The second Reidemeister move is an addition/removal of a bigon,
formed by a pair of edges, which are adjacent (not opposite) to each
other at each of the two vertices, see Fig. \ref{2r}.

\begin{figure}
\centering\includegraphics[width=200pt]{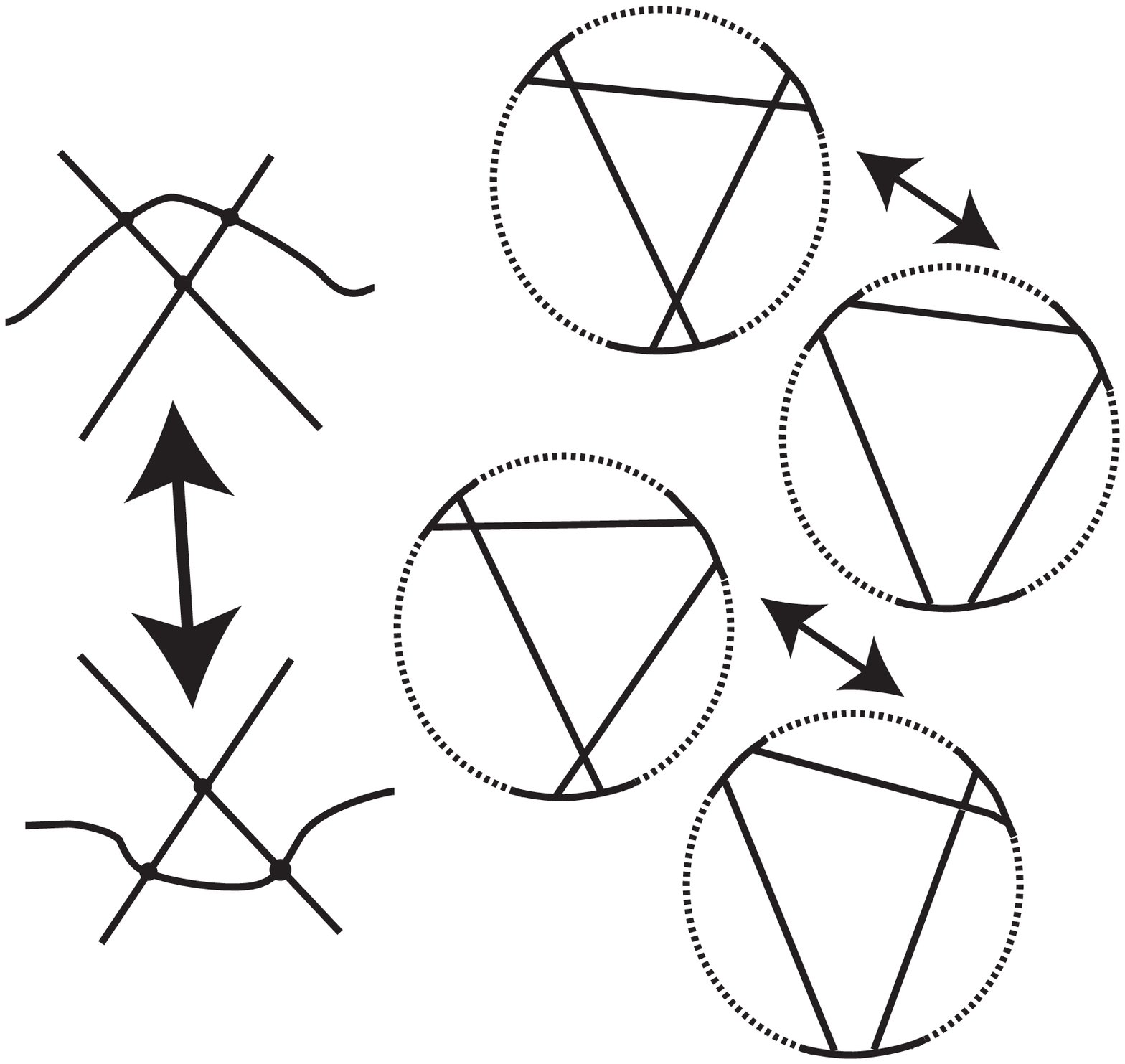} \caption{The Third
Reidemeister Move and Its Chord Diagram Version} \label{3r}
\end{figure}

The third Reidemeister move is shown in Fig. \ref{3r}.

 \begin{dfn}
A {\em free link} is an equivalence class of framed four-valent
graphs modulo Reidemeister moves. Obviously, the number of
components of a framed four-valent graph does not change under
Reidemeister moves, thus, one can talk about the number of
components of a framed links. A {\em free knot} is a $1$-component
free link.

Free knots can be treated as equivalence classes of corresponding
Gauss diagrams (chord diagrams) modulo corresponding moves on Gauss
diagrams (which mimic the moves on four-valent framed graphs).
Analogously, for free-links one can introduce Gauss diagrams on
several circles which stand for different components of the link,
and chord ends may belong to different circles.

Here a chord diagram (Gauss diagram) on one circle is a split
collection of unordered pairs of distinct points on a circle
 $S^{0}_{1}\sqcup \dots \sqcup S^{0}_{k}\subset S^{1}$.
The circle $S^{1}$ is here called the {\em cycle} of the chord
diagram.

Every pair of points $S_{i}^{0}$ is called a {\em chord} and points
from this pair are called chord ends. We call two chords
$S_{i}^{0},S_{j}^{0}$  {\em linked} if the chord ends $S_{i}^{0}$
belong to different connected components of  $S^{1}\backslash
S_{j}^{0}$.

A chord is {\em even} if the number of chords linked with it, is
even, and {\em odd} otherwise.
\end{dfn}

By an {\em even symmetric configuration} $C$ on a chord diagram $D$
we mean a set of pairwise disjoint segments $C_{i}$ on the cycle of
the chord diagram which possess the following properties:

1) the ends of the segments do not coincide with chord ends;

2) the number of chords inside any segment is finite;

3) the number of endpoints of chords inside any segment is even;

4) every chord having one endpoint in $C$ has the other endpoint in
$C$.

5) Consider the involution $i$ of the cycle which fixes all points
outside the segments $C_{i}$ and reflects all segments along the
radii. This involution naturally defines new points to be chord
ends, and thus defines the new chord diagram $i(D)$.

We require that the configuration $C$ is symmetric, i.\,e., chord
diagrams  $D$ and $i(D)$ are equal.

Note that if a chord has two endpoints in different segments $C_{i}$
and $C_{j}$ then it is not self-symmetric.

By an {\em elementary cobordism} we mean a transformation of a chord
diagram deleting all chords belonging to the even symmetric
configuration, as well as the inverse transformation.

We say that two Gauss diagrams are {\em cobordant} if one can be
obtained from the other by a sequence of elementary cobordisms and
third Reidemeister moves.

 \begin{rk}
The first two Reidemeister moves are partial cases of elementary
cobordisms, unlike the third Reidemeister move.
 \end{rk}

Since the first two Reidemeister moves are partial cases of
elementary cobordisms, it makes sense to talk about {\em cobordism
classes of free knots}, not only {\em cobordism classes of Gauss
diagrams}.

 \begin{rk}
The definition of cobordisms given above agrees with the definition
of {\em word cobordism} (nanoword cobordism) \cite{Tu}. With each
(cyclic) double occurrence word in some given alphabet one
associates a Gauss diagram, and words which are cobordant in
Turaev's sense yield cobordant diagrams. However, {\em nanowords}
usually correspond to Gauss diagrams with some {\em decoration
{\rom(}labeling\/{\rom)}} of chords, and the notion of elementary
cobordism usually requires some conditions imposed on chords
belonging to the symmetric configuration.

In this sense, cobordism classes of free knots are the simplest
variant of cobordism classes of nanowords. Nevertheless, we show
that there are free knots which are not cobordant to zero.
 \end{rk}

The main result of the present paper is to prove the existence of
(in fact, infinitely many) free knots which are not cobordant to
zero.

To prove this problem, we shall construct a cobordism invariant for
free knots.

\section{The Mapping $\Delta$ and Its Iterations}

Let $K$ be a framed four-valent graph. Each unicursal component
$K_{i}$ of $K$ can be treated as a four-valent framed graph with a
Gauss diagram $D_{i}$. Thus, some vertices of the graph $K$ can be
represented by chords of one of $D_{i}$'s (namely, those vertices
lying on one unicursal component). Among these, let us choose {\em
even} vertices (in the sense of the Gauss diagram $D_{i}$), and at
each even vertex $X$ of $K$, we consider the smoothing $K_{X}$ (one
of the two possible) for which the number of unicursal components is
greater than that of $K$ by $1$.

Now let ${\cal G}$ be the set of $\Z_{2}$-linear combinations of
equivalence classes of framed four-valent graphs modulo second and
third Reidemeister moves.

Set $\Delta(K)=\sum\limits_{X} K_{X}\in {\cal G}$, where the sum is
taken over all crossings $X$ of $K$.

Then the following statement holds.

 \begin{st}
The map $\Delta$ is a well defined map from ${\cal G}$ to ${\cal
G}$.
 \end{st}

 \begin{proof}
Indeed, assume $K$ is obtained from $K'$ by a third Reidemeister
move. Consider the three crossings $a_{1},a_{2},a_{3}$ of $K$
involved in this move, and the corresponding crossings
$a'_{1},a'_{2},a'_{3}$ of $K'$. By construction, $a_{i}$ lies in one
unicursal component of $K$ if and only if $a'_{i}$ lies in one
unicursal component of $K'$. Moreover, $a$ is even if and only if
$a'_{i}$ is even.





It is now easy to see that whenever $a_{i}$ is even, the smoothing
$K_{a_i}$ gives the same impact to ${\cal G}$ as that of $K'_{a'_i}$
(the corresponding framed 4-valent graphs are either isomorphic or
differ by a second Reidemeister move).

Now, if $K$ and $K'$ differ by a second Reidemeister move and $K'$
has two more crossings $a,b$ in comparison with $K$, then it
obviously follows that either both crossings $a$, $b$ are even or
none of them is even. In the first case, the summands in $\Delta(K)$
are in one-to-one correspondence with those in $\Delta(K')$ and the
corresponding diagrams in each pair differ by a second Reidemeister
move. If both $a$ and $b$ are odd, then it is obvious that the
smoothings at these crossings give equal impact to $K'$, and since
we are working over $\Z_{2}$, they cancel each other.
 \end{proof}

Consequently, if we take $\Delta$ $k$ times, the resulting map
$\Delta^{k}$ is also invariant.

So, if $K$ and $K'$ are two framed four-valent graphs which are
obtained from each other by a third Reidemeister move, then for
every positive integer  $k$ we have $\Delta^{k}(K')=\Delta^{k}(K)$.

\section{The map $\Gamma$}

Let $L$ be a framed four-valent graph with $k$ unicursal components.
With  $L$ we associate a graph $\Gamma(L)$ (not necessarily
four-valent, but without loops and multiple edges) and a number
$I(L)$ according to the following rule.

The graph $\Gamma(L)$ will have $k$ vertices which are in one-to-one
correspondence with unicursal components of  $L$. Two vertices are
connected by an edge if and only if the corresponding components
share an odd number of points.

The following statement is evident.

 \begin{st}
If two framed four-valent graphs $L$ and $L'$ are homotopic then
$\Gamma(L)=\Gamma(L')$.\label{st1}
 \end{st}

Now, we define the number $j(L)$ {\em from the graph $\Gamma(L)$} in
the following way.

If $\Gamma(L)$ is not connected we set $j(L)=0$, otherwise $j(L)$ is
set to be the number of edges of $\Gamma(L)$.

\section{The invariant}

Fix a natural number $n$.

Let $\cal G$ be a linear space generated over $\Z_{2}$ by formal
vectors $\{a_{i}\}, i \in\mathbb{N}$.

Let  $K$ be a framed four-valent graph with one unicursal component,
and let $\Delta^{n}(K)$ be the corresponding linear combination of
four-valent graphs.

For a four-valent framed graph $L$, we set $I(L)=a_{j(L)}$ if
$j(L)\neq 0$ and $I(L)=0$ otherwise. We extend this map to
$\Z_{2}$--linear combinations of framed four-valent graphs by
linearity.

Now, set
 $$
\operatorname{I}^{(n)}(K)=I(\Delta^{n}(K)).
 $$

The main result of the paper is the following

 \begin{theorem}
If  $K$ and $K'$ are cobordant then
$\operatorname{I}^{(n)}(K)=\operatorname{I}^{(n)}(K')$.\label{thmm}
 \end{theorem}

The proof of this theorem follows from two statements. By virtue of
Statement~\ref{st1}, the mapping $\operatorname{I}^{(n)}(\cdot)$ is
invariant with respect to the third Reidemeister move (since so is
$\Delta^{n}$).

Moreover, the following statement holds

 \begin{st}
If $K_{2}$ is obtained from $K_{1}$ by an elementary cobordisms
{\rom (}removal of an even symmetric configuration\/{\rom)}, then
$\operatorname{I}^{(n)}(K_{1})=\operatorname{I}^{(n)}(K_{2})$\label{st3}
 \end{st}

Having proved Statement~\ref{st3}, we shall get Theorem~\ref{thmm}.
Indeed, it will follow that $\operatorname{I}^{(n)}$ is invariant
under all elementary cobordisms and third Reidemeister moves, hence,
under arbitrary cobordisms.

We will prove Statement~\ref{st3} in the last section of our paper

\section{An example}

\begin {figure}
 \centering\includegraphics[width=150pt]{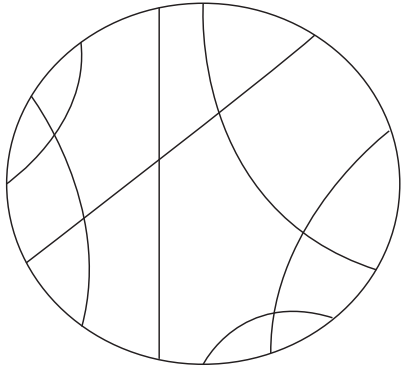}
  \caption{The Gauss diagram of a free knot not cobordant to the unknot}
  \label{chord}
 \end {figure}

  \begin {examp}
Consider the free knot $K$ represented by the Gauss diagram shown in
Fig~\ref {chord}. Let us calculate $\operatorname{I}^{(3)}(K)$. We
have:
  \begin {gather*}
\Delta(\chord)=\chordff+\chordfs+\chordft,\\
\Delta^{2}(\chord)=\Delta(\chordff)+\Delta(\chordfs)+\Delta(\chordft)\\
=\chordsf+\chordss+\chordst+\chordsfo+\chordsfi+\chordssi,\\
\Delta^{3}(\chord)=\Delta(\chordsf)+\Delta(\chordss)+\Delta(\chordst)\\
+\Delta(\chordsfo)+\Delta(\chordsfi)+\Delta(\chordssi)\\
=\chordtf+\chordts+\chordtt+\chordtfo+\chordtfi+\chordtsi\\
=\chordtf+\chordts.
 \end {gather*}
Thus, applying $\Delta^{3}$ we get
$\operatorname{I}^{(3)}(K)=a_{4}$.

Thus, by Theorem~\ref{thmm}, the cobordism class of $K$ is
non-trivial.
  \end {examp}

Obviously, there are infinitely many such examples. We shall
construct them and discuss the {\em cobordism group} of free knots
as well as further invariants in a separate publication

\section{Sketch of the Proof of Statement~\ref{st3}}

Let $K_{2}$ be the Gauss diagram obtained from a Gauss diagram
$K_{1}$ by deleting an even symmetric configuration ${\tilde K}$.

The chords of $K_{1}$ belong to three sets:

1. the set of those chords which corresponding to chords of $K_{2}$,
we denote them for both $K_{1}$ and $K_{2}$ by $\gamma_{j}$'s.

2. the set of those chords $\beta_{j}$ which are fixed under the
involution $i$ on ${\tilde K}$,

\newcommand{\az}{{\bar {\alpha}}}

3. the set of pairs of chords $\alpha_{k}$ and
$\az_{k}=i(\alpha_{k})$ which are obtained from each other by the
involution $i$ (here ${\bar {\bar \alpha}}_{k}=\alpha_{k}$).

Now, for every framed four-valent graph $L$, $\Delta^{n}(L)$ is a
sum of some subsequent smoothings $\Delta^{(p_{1}\dots p_{k})}(L)$,
where all $p_{i}$'s are vertices of $L$ where $p_{i}$ occurs to be
even after smoothing all $p_{1},\dots, p_{i-1}$.

Now, $\Delta^{n}(K_{1})$ naturally splits into $3$ types of
summands:

1. Those where all $p_{i}$ are some $\gamma_{j}$. These smoothings
are in one-to one correspondence with smoothings of $K_{2}$. We
claim that the corresponding elements $I(\Delta^{(p_{1}\dots
p_{k})}(K_{1}))$ and $I(\Delta^{(p_{1}\dots p_{k})}(K_{2}))$ are
equal.

2. Those where at least one of $p_{i}$ is $\beta_{j}$, and neither
$\alpha$'s nor $\az$'s occur among $p_{i}$. We claim that each of
these summands $I(\Delta^{(p_{1}\dots p_{k})}(K_{1}))$ is zero.

3. Those summands where at least one of $p_{i}$'s is $\alpha_{j}$ or
$\az_{j}$.

These summands are naturally paired: the elements
$I(\Delta^{(p_{1}\dots p_{k})}(K_{1}))$ and $I(\Delta^{({\bar
p}_{1}\dots {\bar p}_{k})}(K_{1}))$ are equal.

Let us first prove 1. Since segments of our even symmetric
configurations have no common points with chords $\gamma_{j}$, we
see that after smoothing along any of $\gamma$'s, every segment will
completely belong to one circle. This means that the corresponding
graphs $\Gamma(\Delta^{(p_{1}\dots p_{k})}(K_{1}))$ and
$\Gamma(\Delta^{(p_{1}\dots p_{k})}(K_{2}))$ are isomorphic.

Indeed, the vertices corresponding to $\beta$  do not change the
graph at all, since every vertex corresponding to $\beta$ is an
intersection of one component of $\Delta^{\dots}$ with itself.

Moreover, the vertices $X_{i}$ and $X'_{i}$ corresponding to
$\alpha_{i}$ and $\az_{i}$ belong to the same pair of component, so
their impact to the graph cancels.

To prove 2, let us take one $p_{j}=\beta_{k}$, and consider the
segment $C_{i}$ of $K_1$ where $\beta_{k}$ lies.

Without loss of generality, we may assume that  $\beta_{k}$ is the
innermost chord in $K_1$ amongst those chords $p_{l}$ we use for
smoothings.

Now, our summand looks like $\Delta^{\dots \beta_{k}\dots}$.
Smoothing $K_{1}$ along $\beta_{k}$ cuts a free knot (component)
which has ends only in the segment $C_{i}$. It is obvious that this
unicursal component will be split in the sense of the graph
$\Gamma$: it will have even intersection with any other unicursal
component (since whenever it shares a vertex $X$ of some other
component, this vertex necessarily corresponds either to some
$\alpha$ (or to some $\az$) and the corresponding $\az$ (resp.,
$\alpha$) will be shared by the same two components).

The proof of 3 follows from one basic fact: {\em the number of
components of a $1$-manifold obtained from a chord diagram by
smoothing some of its chords can be defined from the intersection
graph of this chord diagram}, see~\cite{Sob}. We shall give a
rigorous proof of 3 in a separate publication.

\vspace{1cm}

The authors express their gratitude to V.P.Ilyutko for implementing
a computer program calculating the invariant $\operatorname{I}$.



\begin{thebibliography}{100}

\bibitem[Gib]{Gib} Gibson, A., {\em Homotopy Invariants of Gauss
Words}, ArXiv:Math.GT/0902.0062.

 \bibitem[Ka]{Kau}
L.\,H.~Kauffman, Virtual knot theory, Eur. J. Combinatorics. 1999.
V.\, 20, N.\, 7, pp.\, 662--690.

\bibitem[Ma]{Ma} Manturov, V.O., {\em On Free Knots},
ArXiv:Math.GT/0901.2214

 \bibitem[Sob]{Sob}
E.~Soboleva, Vassiliev Knot Invariants Coming from Lie Algebras and
$4$-Invariants (2001),  {\em Journal of Knot Theory and Its
Ramifications}, {\bf 10} (1), pp.\ 161--169.

\bibitem[Tu]{Tu} Turaev, V.G., {\em Cobordisms of Words},
Arxiv:Math.CO/0511513, v.2.

\end{thebibliography}
 \end {document}